\documentclass[12pt]{article}
\date{}
\title{One-dimensional symmetry for solutions of  Allen Cahn fully nonlinear equations.}

\author{ F. Demengel\, I. Birindelli}

\newtheorem{theo}{Theorem}[section]
\newtheorem{prop}[theo]{Proposition}
\newtheorem{rema}[theo]{Remark}
\newtheorem{defi}[theo]{Definition}

\newtheorem{exe}{Example}[section]

\newcommand{\ph}{\varphi}
\newcommand{\R}{{\rm I}\!{\rm  R}}
\def\R{{\rm I}\!{\rm  R}}

\def\la1{\lambda_1}
\def\grad{\nabla}
\def\ph1{\varphi_1}
\newcommand{\N}{{\bf N}}

\begin{document}
\maketitle
\begin{abstract}
 This article presents some qualitative results for solutions of the fully nonlinear elliptic equation 
 $F(\nabla u, D^2 u) + f(u)=0$ in $\R^N$.  Precisely  under some additional assumptions  on $f$, if    $-1\leq u\leq 1$ and $\lim _{x_1\rightarrow \pm \infty}  u(x_1, x^\prime) = \pm 1$ uniformly with respect to $x^\prime$, then the solution depends only on $x_1$.   \end{abstract}
 \section{Introduction}
 The sliding method was crystalized in \cite{BN}  by Berestycki and Nirenberg in  order to prove monotonicity of solutions of
 \begin{equation}\label{eqq1} \Delta u +f(u)=0\quad \mbox{ in } \Omega\subset\R^N.
 \end{equation}
This powerful method uses two features  of the Laplacian, comparison principle and invariance with respect to translation. The idea is: Fix a direction;
first slide in that direction enough for the intersection of the slided  domain with $\Omega$ to be small enough or "narrow enough" for the maximum principle to hold. This  allows  to compare the value of the solution at different points of the domain.  Then continue "sliding" until reaching a critical position.

Coupling simplicity with ductility, the sliding method of \cite{BN} has been incredibly influential, it is possible to count over two hundred citations of the work (e.g. through google scholar). We shall here only recall the work by Berestycki, Hamel and Monneau \cite{BHM} where the technic is  used to prove the so called Gibbons conjecture . This was simultaneously and independently solved by Barlow, Bass and Gui \cite{BBG} and Farina   \cite{F}. Precisely in \cite{BHM}, they prove that if $f$ is a $C^1([-1,1])$ function decreasing near $-1$ and $1$, with  $f(-1)=f(1)=0$ (typically, $f(u)=u-u^3$) then the solutions of
(\ref{eqq1}) in $\R^N$ that converge uniformly to 1 or -1 at infinity in some fixed direction, say $x_1$, are in fact one dimensional i.e functions of $x_1$ alone.
In \cite{BHM}, the sliding method is coupled with a maximum principle (comparison principle) in unbounded domains contained in some cone.

As is well known the Gibbons conjecture is a weak form of the famous  De Giorgi's conjecture which states that for $f(u)=u-u^3$, the level sets of monotone, entire solutions of (\ref{eqq1}) are hyperplanes for $N\leq 8$.  This result has been proved in dimension 2 and 3 respectively by Ghoussoub and Gui \cite{GG} and by Ambrosio, Cabr\'e \cite{AC}, while Del Pino, Kowalcyk and Wei \cite{DKW} have proved that it does not hold for $N>8$ by constructing a counter example. 
Savin has proved the case $4\leq N\leq 8$, with the further condition that
the limit be $\pm 1$ in a direction at infinity, in that case this condition is not assumed to be uniform with respect to the other variables. See also \cite{VSS} for analogous results concerning the $p$-Laplacian.

In the present note we extend Gibbons conjecture to fully nonlinear operators. Precisely, we consider entire bounded solutions of 
   \begin{equation}\label{eq1}
  F(\nabla u, D^2 u) + f(u)=0\quad  \mbox{in}\quad \R^N,
  \end{equation}
where $F(\nabla u, D^2 u):=|\grad u|^\alpha \tilde F(D^2u)$ with $\alpha>-1$ and $\tilde F$ is uniformly elliptic.
With the same conditions on the nonlinearities of $f$ as in \cite{BHM}, we prove that for any solution such that
$\lim _{x_1\rightarrow \pm \infty}  u(x_1, x^\prime) = \pm 1$ uniformly with respect to $x^\prime$ and such that $|\grad u |>0$ in $R^N$ then $\partial_{x_1} u\geq 0$ and 
$u$ is a function of $x_1$ alone.

Many remarks are in order. Let us note that in the case $\alpha\leq 0$, some recent regularity results \cite{BD9} prove that locally Lipschitz solutions are in fact  ${\cal C}^{1, \beta}$ for some $\beta <1$,  and this regularity is sufficient  to prove the results enclosed here. For $\alpha>0$ the  ${\cal C}^{1}$ regularity is a consequence of the hypothesis on the positivity of the norm of the gradient.

A key ingredient in the proof of this result, which is of  independent interest, is the following, strong comparison principle.
  \begin{prop}\label{strict0}
  Suppose that $\Omega$ is some open  set, and $x_o, r$ such that $B(x_o, r) \subset \Omega$. 

   Suppose that $f$ is ${\cal C}^1$ on $\R$ , and that $u$ and $v$ are, respectively,  ${\cal C}^1$  bounded  sub- and super-solutions of 
   $$F(\nabla w, D^2 w)+ f(w) =0\quad \mbox{in}\quad \Omega$$
such that $u\geq v$ and $\nabla v \neq 0$ (or $\nabla  u\neq 0$) in $B(x_o, r)$,  then, 
  either  $u> v$  or $u\equiv v$ in $B(x_o, r)$. 
  \end{prop}
 Observe that the condition  that the gradient needs to be different from zero cannot be removed. Indeed, for any $m,k\in Z$ with $k\leq m$ the functions
 $$u_{k,m}(x)=\left\{\begin{array}{lc}
 1 & \mbox{for}\ x_1\geq (2m+2)\pi\\
 \cos x_1  & \mbox{for}\ (2k+1)\pi\leq x_1\leq (2m+2)\pi\\
   -1 & \mbox{for}\ x_1\leq (2k+1)\pi
   \end{array}
   \right.
   $$
 are viscosity  solutions of 
 $$|\grad u|^2(\Delta u)+(u-u^3)=0,$$ and they are ${\cal C}^{1,\beta}$ for all $\beta <1$. 
 
Observe that e.g. $u_{0,0}\geq u_{0,i}$ for all $i\geq 1$  and $u_{0,0}(2\pi,y)=u_{0,i}(2\pi,y)$ but the functions don't coincide.

This example suggests that there may be solutions that are not one dimensional if the condition on the gradient is  removed.

\bigskip
When $\alpha=0$, De Silva and Savin in \cite{DSS}, have proved the analogue of De Giorgi's conjecture for uniformly elliptic operators in dimension 2. With $f$ as above, they prove that if there exists a one dimensional monotone solution i.e. $g:\R\mapsto [-1.1]$ such that $u(x)=g(\eta\cdot x)$ is a solution of
\begin{equation}\label{eqq2}
\tilde F(D^2u)+f(u)=0\quad\mbox{in}\quad\R^2
\end{equation}
satisfying $\lim_{t\rightarrow\pm\infty}g(t)=\pm1$
then, all monotone bounded solutions of (\ref{eqq2}) are one dimensional, i.e. their level sets are straight lines.

Let us mention that without any further assumptions on $f$ solutions may not exists. Indeed, let $\tilde F(D^2u)= {\cal M}^+_{a,A} (D^2 u)$ where for any symmetric matrix $M$ with eigenvalues $e_i$,
$${\cal M}^+_{a,A}(M)=a\sum_{e_i<0} e_i +A\sum_{e_i>0} e_i.$$
Then, as shown in the last section, for $a<A$ there are no one dimensional solutions of
$$  {\cal M}^+_{a,A} (D^2 u)+u-u^3=0,$$
that satisfy the asymptotic conditions.
In that section we study conditions on $f$ that guarantee existence of solutions of the ODE 
$$  |u^\prime |^\alpha {\cal M}^+_{a,A} (u^{\prime\prime})+f(u)=0$$
that satisfy $\lim_{x\rightarrow\pm\infty}u(x)=\pm 1$.

\bigskip

While completing this work, we have received a paper by Farina and Valdinoci,  \cite{FV}, who treats Gibbons conjecture in a very general setting that includes the case $\alpha = 0$.

 \section{ Assumptions and known results}

In the whole paper we shall suppose the following hypotheses on the operator $F$. 

Let $S$ be the set of $N\times N$ symmetric matrices, and let  $\alpha >-1$. Then $F$  is defined on $\R^N\setminus\{0\}\times S $  by
\begin{equation}\label{deff}
F( p, M) = |p|^\alpha  \tilde F( M),
\end{equation}
where $\tilde F$ satisfies

$$\tilde F(tM)=t\tilde F(M)\quad\mbox{ for any }\quad t\in \R^+, M\in S,$$
and there exist $A\geq a>0$ such that for any $M$  and any $N\in S$ such that $N\geq 0$ 
\begin{equation}\label{eqaAF}
a tr(N)\leq \tilde F(M+N)-\tilde F(M) \leq A tr(N).
\end{equation}

 \begin{exe}
   
   1) Let  $0< a < A$ and ${\cal M}_{a,A}^+(M)$ be the Pucci's operator 
   ${\cal M}_{a,A}^+ (M) = Atr(M^+)-a tr(M^-)$ where $M^\pm$ are the positive and negative part    of $M$, and 
    ${\cal M}_{a,A}^-(M)=-  {\cal M}_{a,A}^+ (-M)$. 
       Then $F$ defined as  
   
   $$F( p,M) = |p|^\alpha {\cal M}_{a,A}^\pm (M) $$
 satisfies the assumptions.
   
2) Let $B$ be a symmetric positive definite  matrix then     $F( p,M) = |p|^\alpha(tr(BM) )$, is another example of operator satisfying the assumptions.

    \end{exe} 
We now recall what we mean by viscosity solutions in our context : 
  
     \begin{defi}\label{def1}

 Let $\Omega$ be a bounded domain in
$\R^N$, let $g$ be a continuous function on $\Omega\times \R$, then
$v$,   continuous on $\overline{\Omega}$ is called a viscosity super-solution (respectively sub-solution) 
of
$F(\grad u,D^2u)=g(x,u)$ if for all $x_0\in \Omega$, 

-Either there exists an open ball $B(x_0,\delta)$, $\delta>0$  in $\Omega$
on which 
$v$ is  a constant $ c
$ and 
$0\leq g(x,c)$, for all $x\in B(x_0,\delta)$
(respectively 
$0\geq g(x, c)$ for all $x\in B(x_0,\delta)$) 

-Or
 $\forall \varphi\in {\mathcal C}^2(\Omega)$, such that
$v-\varphi$ has a local minimum (respectively local maximum) at  $x_0$ and $\grad\varphi(x_0)\neq
0$, one has
$$
F( \grad\varphi(x_0),
 D^2\varphi(x_0))\leq g(x_0,v(x_0)).
$$
(respectively  
$$F( \grad\varphi(x_0),
 D^2\varphi(x_0))\geq g(x_0,v(x_0))).$$

A	 viscosity solution is a function which is both a super-solution and a sub-solution.
\end{defi}

\begin{rema}
 When $F$ is continuous in $p$,  and $F(0,0)=0$, this  definition is equivalent to the classical definition
of viscosity solutions, as in the User's guide \cite{CIL}.
\end{rema}
We now give a definition that will be needed in the statement of our main theorem.

\begin{defi}\label{ddd} We shall say that $|\grad u|\geq m>0$ in $\Omega$ in the viscosity sense, if for  all  $ \varphi\in {\mathcal C}^2(\Omega)$, such that
$u-\varphi$ has a local minimum or a local maximum at some  $x_0\in\Omega$,
 $$|\grad\varphi(x_0)|\geq m.$$
\end{defi}

 In our context, since the  solutions considered have their gradient different from zero everywhere, the  viscosity solutions  can be intended in  the classical meaning.

 We begin to recall some of the results obtained in \cite{BD5} which will be needed in this article.

 \begin{theo}\label{thcomp1}
 Suppose that   $c$  is a continuous and bounded function satisfying  $c\leq 0$. 

Suppose that $f_1$ and $f_2$ are continuous and bounded and that 
$u$ and $v$  satisfy
\begin{eqnarray*}
 F( \nabla u, D^2 u)+c(x)|u|^\alpha u& \geq & f_1\quad 
\mbox{in}\quad \Omega, \\ 
F( \grad v,D^2 v)
+ c(x) |v|^{\alpha}v & \leq & f_2 \quad  \mbox{in}\quad
\Omega , \\ 
u \leq  v &&   \quad  \mbox{on}\quad \partial\Omega.
\end{eqnarray*}
 If $f_2< f_1$ then $u \leq v$ in $\Omega$. 
Furthermore, if $c<0$ in $\Omega$  and $f_2\leq f_1$ then $u \leq v$ in $\Omega$. 
\end{theo}

   \begin{prop}\label{remhopf}
 Suppose that  ${\cal O}$ is a smooth bounded domain. Let 
 $u$ be   a solution of 
    \begin{equation}\label{123}
  F( \nabla u, D^2 u) \leq 0 
\quad \mbox{in }\quad {\cal O}. \end{equation}
If  there exists some constant $c_o$,  such that $u>c_o$ inside ${\cal O}$ and $u(\bar x)=c_o$ with $\bar x\in \partial {\cal O}$,  then
$$\liminf_{t\rightarrow 0^+} \frac{u(\bar x-t\vec n)-u(\bar x)}{t}>0,$$ 
where $\vec n$ is the outer normal to $\partial {\cal O}$ at $\bar x$. 
\end{prop}

\begin{rema} \label{remrem} In particular Proposition \ref{remhopf} implies that a non constant  super-solution of (\ref{123}) in a domain $\Omega$ has no interior minimum.

 If $c_o = 0$, the result can be extended  in the following manner : Suppose that $\beta\geq \alpha $,    that  $c$ is continuous and bounded, and $u$ is a nonnegative solution of 
 $$F(\nabla u, D^2 u) + c(x) u^{1+\beta}  \leq 0$$
 then  either $u\equiv 0$ or $u>0$ in $\Omega$. In that last case, if $u = 0$ on some point $x_o\in \partial \Omega$, then $\partial_{\vec n} u(x_o)>0$. 
 \end{rema}

   We now recall the regularity results obtained in \cite{BD9}. 

 \begin{theo}\label{reg}
 
 Suppose that  $\Omega$ is  a bounded ${\cal C}^2$ domain  and 
 $\alpha \leq 0$. Suppose that  $ g$ is continuous on $\Omega \times \R$ . Then the  bounded solutions  of 
  \begin{equation} \label{012}\left\{ \begin{array}{lc}
  F(\nabla u, D^2 u) = g(x, u(x)) &\mbox{ in }\    \Omega,\\
  u = 0 &\ \mbox{ on } \ \partial \Omega,
  \end{array}\right. 
  \end{equation}
  satisfy
  $u\in {\cal C}^{1,\beta}(\overline{\Omega})$, for some $\beta \in (0, 1)$ . 

Furthermore if $\Omega$ is a domain (possibly unbounded) of  $\R^N$  and if $u$ is bounded and locally Lipschitz then $u\in {\cal C}_{loc}^{1,\beta}(\Omega)$ for some $\beta\in (0,1)$.
  \end{theo}

  When $\alpha>0$, ${\cal C}^1$ regularity results are not known except for the one dimensional case or the radial case,  however here, since the solutions that we consider have the gradient bounded away from zero, this regularity is  just a consequence of classical results and a priori estimates.
  Indeed next theorem is just an application of Theorem 1.2 of \cite{CDV}, which in turn is the extension of Caffarelli's classical result:
\begin{theo}\label{IA}
Suppose that $\Omega$  is a  (possibly unbounded) domain, and  that  $g$ is  ${\cal C}^1$ and  bounded.  Let   $u$ be  a  bounded solution of 
\begin{equation}\label{1a}
   F(\nabla u, D^2 u) = g(u)\ \mbox{ in }\    \Omega.
   \end{equation}
If  $|\grad u|\geq m>0$ in $\Omega$  in the sense of Definition \ref{ddd},   there exists  $\beta\in (0,1)$ and $C=C(a,A,N,|g(u)|_\infty, m)$ such that
if $B(y, \rho) \subset \Omega$, 
    \begin{equation} \label{eqholdloc} \|u\|_{{\cal C}^{1,\beta}(B(y, \frac{\rho}{2}))}\leq C \sup_ {B(y, \rho)}|u|.
    \end{equation}

\end{theo}
{\em Proof.}  We introduce the operator:
$$G(v,\grad v, D^2v):=\tilde F(D^2v)-g(v)\sup\left(|\grad v|,\frac{m}{2}\right)^{-\alpha}.$$
If $u$ is a solution of (\ref{1a}) such that in the viscosity sense $|\grad u|\geq m>0$,  then it is a solution of 
$$ G(u,\grad u, D^2 u)=0\quad \mbox{in}\quad \Omega.$$
Indeed, e.g. if $\varphi\in {\cal C}^2$ is such that 
$(u-\varphi)(x)\geq (u-\varphi)(\bar x)$  for some $\bar x\in\Omega$, then $|\nabla \varphi |(\bar x)\geq m$ and
$$|\nabla \varphi|^\alpha(\bar x) \tilde F(D^2\varphi(\bar x))\geq g(u(\bar x))\Rightarrow \tilde F(D^2\varphi (\bar x))- |\nabla \varphi (\bar x)|^{-\alpha} g(u(\bar x))\geq 0.$$
In order to apply Theorem 1.2 of \cite{CDV}, it is enough to remark that,  $G$ does not depend on $x$ and therefore the condition on the modulus of continuity is automatically satisfied. 

Furthermore,
 the dependence on the gradient is Lipschitz, where the Lipschitz constant depends on $m$ and  $|g(u)|_\infty$. 
 Applying Theorem 1.2 of \cite{CDV} we have obtained the above estimate and
$u\in C^{1,\beta}(\Omega)$. This ends the proof.
      
 \section{Comparison principles}

As mentioned in the introduction, we begin by proving a strong comparison principle, that extends the one obtained in \cite{BD9}.
 
  \begin{prop}\label{strict}
  Suppose that $\Omega$ is some open subset of $\R^N$,  
$f$ is ${\cal C}^1$ on $\R$ . Let $u$ and $v$ be  ${\cal C}^1$  bounded  sub-solution  and super-solution of 
   $$F(\nabla u, D^2 u)+ f(u) =0\quad \mbox{in}\ \Omega.$$
Suppose that  ${\cal O}$ is  some connected subset of $\Omega$, with   $u\geq v$ and $\nabla v \neq 0$ (or $\nabla  u\neq 0$) on  ${\cal O}$ ,  then 
  either  $u> v$  or $u\equiv v$ in ${\cal O}$. 
\end{prop}
\begin{rema} Of course when $\alpha=0$ the strong comparison principle is classical and holds without requiring that the gradient be different from zero.
\end{rema}
\noindent{\em Proof of Proposition \ref{strict}.}
 We write the proof in the case $\alpha< 0$, the  changes  to bring when $\alpha>0$  being obvious.

  We argue as in \cite {BD9}. Suppose that $x_o$ is some point where $u(x_o)>v(x_o)$ (if such point doesn't exist we have nothing to prove). 
  
 Suppose by contradiction that there exists some point $x_1$ such that $u(x_1) = v(x_1)$. It is clear that it  can be chosen in such a way that, for  $R = |x_1-x_o|$,
 $u>v$ in $B(x_o, R)$ and $x_1$ is  the only point in the closure of that ball  on which $u$ and $v$ coincide.    
Without loss of generality, one can assume that $B(x_o, {3R\over 2}) \subset {\cal O}$.  
  
  We can assume without loss of generality that $v$ is the function whose gradient is bounded away from zero. 
 Let then  $L_1 = \inf_{B(x_o, {3R\over 2})} |\nabla v|>0$, $L_2 = \sup_{B(x_o, {3R\over 2})} |\nabla v|$.     
  We  will  prove that there exist  two constants  $c>0$ and   $\delta >0  $ such that 
  $$ u \geq v + \delta  ( e^{-c|x-x_o|}- e^{-3cR\over 2})\equiv v+ w\quad \mbox{in}\quad {R\over 2}\leq |x-x_o| = r\leq {3R\over 2}.$$  
This will  contradict  the fact that $u(x_1) = v(x_1)$.

 Let $\delta \leq \displaystyle \min_{|x-x_o|= {R\over 2} }(u-v)$,   so that 
 $$u\geq v+w\quad  \mbox{on} \quad
 \partial\left(B(x_o, {3R\over 2})\setminus \overline{B(x_o, {R\over 2}})\right).$$
Define 
$$\gamma(x) = \left\{ \begin{array}{lc}
 {f(u(x))-f(v(x))\over u(x)-v(x)} &\mbox{ if } \ u(x) \neq v(x)\\
  f^\prime (u(x)) & \mbox{ if } \ u(x) = v(x). 
  \end{array}\right.$$
  Since $f$ is ${\cal C}^1$ and the  functions $u$ and $v$ are bounded, $\gamma$ is continuous and bounded. 
  We  write 
 $$
f(u)  =\gamma (x)  (u-v) + f(v), $$
  $$ F(\nabla u, D^2 u) -(|\gamma |_\infty+1)  (u-v) = -f(v)+ (-\gamma-|\gamma|_\infty-1)  (u-v)\leq F(\nabla v, D^2 v). $$
We shall prove that,  for $c$ chosen conveniently, 
   $$F(\nabla v, D^2 v) < F(\nabla (v+ w), D^2(v+ w)) -(|\gamma|_\infty+1) 
 w,$$
    this will imply that 
    $$F(\nabla u, D^2 u)-(|\gamma|_\infty+1)  u \leq F(\nabla (v+ w), D^2(v+ w)) -(|\gamma|_\infty+1)  (v+ w).$$
Let $\varphi$ be some test function for $v$ from above,
a simple calculation on $w$ implies that, if $c \geq {1\over a}({2(2A(N-1) \over R} )$ then
                
     \begin{eqnarray*} 
     |\nabla \varphi+ \nabla w|^\alpha &\cdot&\tilde F (x, D^2 \varphi+ D^2 w)-(|\gamma|_\infty+1)    w\\
     &\geq&    |\nabla \varphi+ \nabla w|^\alpha \tilde F (x,D^2 \varphi)    +   |\nabla \varphi+ \nabla w|^\alpha {\cal M}^- (D^2 w)-(|\gamma|_\infty+1)   w \\
        &\geq& |\nabla \varphi+ \nabla w|^\alpha { F(\nabla \varphi, D^2\varphi)\over |\nabla \varphi|^\alpha} +\\
       && + |\nabla \varphi+ \nabla w|^\alpha\frac{ac^2}{2}  \delta e^{-cr} -(|\gamma|_\infty+1)  \delta e^{-cr}.
\end{eqnarray*}
We also impose  $\delta < {R L_1e\over 16 }$ so that  $|\nabla w|\leq {|\nabla \varphi]\over 8}$; then
the inequalities 

    $$||\nabla \varphi+\nabla w|^\alpha-|\nabla \varphi|^\alpha |\leq |\alpha | |\nabla w|| \nabla \varphi|^{\alpha-1}\left( {1\over 2}\right)^{\alpha-1}\leq {|\nabla \varphi|^\alpha \over 2}$$
imply that  
  $$
     |\nabla \varphi+ \nabla w|^\alpha \left(\tilde F(x,D^2 \varphi+ D^2 w)\right) 
\geq    -f(v) -| f(v)|_\infty|\grad\varphi|^{-1}|\alpha| 2^{1-\alpha}c\delta e^{-cr}   +L_2^\alpha {ac^2\over 4} \delta e^{-cr}.
$$
It is now enough to choose  
$$c\geq  {4A(N-1) \over R} + {|\alpha | |f(v)|_\infty 2^{2-\alpha} \over a L_2^{1+\alpha}}+ \left({16(|\gamma |_\infty+1)\over a L_2^\alpha} \right)^{1\over 2}$$ 
to finally obtain
$$
     |\nabla \varphi+ \nabla w|^\alpha \tilde F(x,D^2 \varphi+ D^2 w)-(|\gamma|_\infty+1  ) w\geq    f(v)+{a  c^2\delta L_2^\alpha e^{-cr}\over 8}-(|\gamma|_\infty+1)  \delta e^{-cr}
$$
i.e.
  $$F(x, \nabla (v+ w), D^2 (v+ w)) -(|\gamma|_\infty+1)   w> F(x,\nabla v,D^2 v).$$
  Hence the comparison principle, Theorem \ref{thcomp1}, gives that 
   $$ u\geq v+w\quad\mbox{in} \ B(x_o,\frac{3R}{2})\setminus \overline{ B(x_o,\frac{R}{2})},$$
the desired contradiction. This ends the proof of Proposition \ref{strict}.
\bigskip

From now $f$ will denote a ${\cal C}^{1}$ function defined on $[-1,1]$, such that $f(-1) =  f(1) =0$,  
 and nonincreasing  on the set  $[-1, -1+\delta]\cup [1-\delta, 1]$ for some $\delta\in ]0,1[$.

Next is a comparison principle in unbounded domains that are "strip" like.

 \begin{prop}\label{propcomp}
 Suppose that $u$ and $v$  are ${\cal C}^1$, have values in $[-1,1]$ and are respectively sub and super solutions of 
  $$F(\nabla w,D^2 w) + f(w) =  0\ \mbox{ in }  \ \R^N$$
 with   $F(\nabla u,D^2 u) \in L^\infty$, $F(\nabla v,D^2 v)\in L^\infty$. 
If  $b, c\in \R$ are such that $b < c$, $\Omega = [b,c]\times \R^{N-1}$, 
   $|\nabla u|$ and $|\nabla v|\geq m>0$   and 
   either $u\leq -1+\delta$ or $v\geq 1-\delta$ in $\Omega$, 
then   
   $$u-v\leq \sup_{\partial \Omega } (u-v)^+.$$
   \end{prop}
Proof of Proposition \ref{propcomp}.  
    
Without loss of generality $f$ can be extended outside of $[-1,1]$ in order that $f$ be still ${\cal C}^1$ , bounded, and nonincreasing after $1-\delta $ and before $-1+ \delta$.  Suppose, to fix the ideas, that $v\geq  1-\delta$ in $\Omega$.

We can also assume that $u\leq v$ on $\partial \Omega$. Indeed, since 
$f$ is decreasing after $1-\delta$,  $w= v+ \sup_{\partial \Omega } (u-v)^+$ is a super-solution which satisfies $F(\nabla w ,D^2 w) \in L^\infty$. Suppose by contradiction that  $\sup_\Omega  (u-v) =   \lambda$
     for some $\lambda >0$. 

By definition of the supremum,  there exists some sequence $(x^k)_k$ such that $ (u-v) (x^k)\rightarrow \lambda$. Eventually extracting from $(x^k)_k$ a subsequence, still denoted $(x^k)_k$, we have $x_1^k\rightarrow \bar x_1\in [b,c]$. For any $x=(x_1,x^\prime)$ let
            $$u^k(x_1, x^\prime) = u(x_1, x^\prime + (x^\prime)^k)$$
       and  
         $$v^k(x_1, x^\prime) = v(x_1, x^\prime + (x^\prime)^k).$$

          By the uniform  estimates \ref{eqholdloc} in Theorem \ref{IA} one can extract from $(u^k)_k$ and $(v^k)_k$ some  subsequences,  denoted in the same way,  such that $u^k\rightarrow \bar u$ and  $v^k\rightarrow \bar v$ uniformly on every  compact set of $[b,c ]\times \R^{N-1}$  and $\bar u$ and $\bar v+ \lambda$ are solutions of 
          $$F(\nabla \bar u, D^2 \bar u)  \geq -f(\bar u),$$
          $$F(\nabla(\bar v+\lambda), D^2 (\bar v+\lambda)) \leq -f(\bar v) \leq -f(\bar v+\lambda).$$
Furthermore,  $\bar u\leq \bar v+ \lambda$, and 
through the uniform convergence on the compact set $[b, c] \times \{0\}^{N-1}$, $\lim_k u^k(\bar x_1, 0) = \lim_k u^k (x_1^k, 0)$ and $\lim_k v^k(\bar x_1, 0) = \lim_k v^k (x_1^k, 0)$.  This implies that 

           \begin{eqnarray*}
                      \bar u(\bar x_1, 0) &=& \lim_k u(x_1^k, 0+ {x^\prime}^k)\\&=&
            \lim_k v( x_1^k, 0+ {x^\prime}^k)+ \lambda =\bar v(\bar x_1, 0)+ \lambda. 
            \end{eqnarray*}
            Now using the fact that $|\nabla u|> m$ and  $|\nabla v|> m$ on $[b,c]\times  \R^{N-1}$,  by passing to the limit one gets that $|\nabla \bar u|\geq m>0$ and $|\nabla \bar  v| \geq m$ on that strip, and  the strong comparison principle in  Proposition \ref{strict},   implies that $\bar u \equiv \bar v+ \lambda$.
             
On the other hand,
$$u(b, x^\prime +  {x^\prime}^k)\leq v(b, x^\prime +  {x^\prime}^k)$$ implies,  by passing to the limit that 
$$\bar u(b, x^\prime )\leq \bar v(b, x^\prime )$$
a contradiction.

 \section{Proof of the one dimensionality.}

   We now state precisely and prove  the main result of this paper: 
  
        \begin{theo}\label{th1}
        Let $f$ be defined on $[-1,1]$, ${\cal C}^1$ and such that $f$ is nonincreasing near $-1$ and $1$, with $f(-1)  = f(1) = 0$. 
Let $u$ be a   viscosity solution of 
          $$F(\nabla u, D^2 u) + f(u) =  0\ \mbox{ in }  \  \R^N,$$
with values in $[-1,1]$.            
 Suppose that  
 $\displaystyle\lim_{x_1\rightarrow 
   \pm \infty} u(x_1, x^\prime) = \pm 1$,   uniformly with respect to $x^\prime $,  
and, if $\alpha\neq 0$, suppose that for any $b<c$ there exists $m>0$ such that 
$|\nabla   u(x)|\geq m>0$ in $[b, c]\times \R^{N-1}$ in the viscosity sense. 

\noindent  Then  $u$ does not depend on $x^\prime$ i.e. $u(x_1, x^\prime ) = v(x_1)$ where 
\begin{equation}\label{dim1wholespace}
 \left\{ \begin{array}{lc}
F( v^\prime e_1 , v^{\prime\prime}e_1\otimes e_1) + f(v) = 0& \ \mbox{ in }  \ \R,\\
   |v|\leq 1, \ \displaystyle \lim_{x\rightarrow \pm \infty} v = \pm1 & \end{array} \right.
   \end{equation}
   and $v$ is increasing.
\end{theo}
{\em Proof of Theorem \ref{th1}.} We proceed analogously to the proof given in \cite{BHM}.
First observe that by Theorem  \ref{IA} the solution $u$ is in ${\cal C}_{loc}^{1,\beta}(\R^N)$, so that the condition on the gradient is pointwise and not  only in the viscosity sense.

Let $\delta $ be such that $f$ is nonincreasing on $[-1, -1+\delta]\cup [1-\delta, 1]$. Define 

 $$\Sigma_M^+:=\{x=(x_1,x')\in\R^N,\  x_1\geq M\}\quad\mbox{
and}\quad \Sigma_M^-:=\{x=(x_1,x')\in\R^N,\ x_1\leq M\}.$$ 
By the  uniform behavior  of the solution in the $x_1$ direction,  there exists $M_1>0$ such that 
$$u(x) \geq 1-\delta \quad \mbox{in}\quad \Sigma^+_{M_1},\quad u(x)\leq -1+\delta  \quad \mbox{in}\quad \Sigma^-_{(-M_1)}.$$

Fix any   $\nu=(\nu_1,\dots,\nu_n)$ such that $\nu_1>0$ and let
$u_t (x) := u(x+t\vec \nu)$.

\noindent {\bf Claim 1} : For $t$ large enough, $u_t\geq u$ in $\R^N$. 

For $x\in\Sigma^+_{(-M_1)}$ and for $t$ large enough, say $t > {2M_1 \over \nu_1}$,
   $$ u(x+t \vec\nu ) \geq 1-\delta\quad\mbox{ and }\quad u^t \geq u\quad\mbox{on}\ x_1 = -M_1.
$$
We begin  to prove   that $u_t\geq u$ in $\Sigma^+_{(-M_1)}$.

\noindent   Suppose by contradiction that $\sup_{\Sigma^+_{(-M_1)}}(u-u_t)=m_o>0$. 

  \noindent  Observe that since $\displaystyle \lim_{x_1\rightarrow +\infty} u = \lim_{x_1\rightarrow +\infty} u_t= 1$ uniformly, there exists $M_2$ such that for 
$x_1> M_2\geq  -M_1$, $|u_t-u|< {m_o\over 2}$. Then   $\sup_{\Sigma^+_{(-M_1)}}(u-u_t)=m_o$ is achieved inside 
$[-M_1, M_2] \times \R^{N-1}$. 

On that strip, by hypothesis, there exists  $m>0$ such that $|\nabla u|, |\nabla u_t |\geq m$, and also $u_t\geq 1-\delta$.  Then one can apply the strong comparison principle  in  Proposition \ref{propcomp}  with $b = -M_1$ and $c = M_2$  and obtain that 
$$u-u_t \leq \sup_{\{x_1 = -M_1\}\cup \{ x_1 = M_2\}}  (u-u_t)^+< {m_o\over 2},$$ a contradiction. Finally we have $u\leq u_t$ in $\Sigma^+_{(-M_1)}$.

\noindent We can do the same in $\Sigma^-_{\{-M_1\}}$ by observing that, in that case, $u\leq -1+\delta$. 

This ends the proof of Claim 1.

 \bigskip

 Let $ \tau = \inf\{ t>0,\ \mbox{such that} \ u_t\geq u\in \R^N\}$, by Claim 1, $\tau$  is finite.

\noindent {\bf Claim 2:} $\tau=0$.

 To prove this claim, we argue by contradiction, assuming  that it is positive.   

\noindent We suppose first that
    $$\eta := \inf _{ [-M_1, M_1] \times \R^{N-1}} (u_\tau-u)>0,$$ and we  prove then that there exists $\epsilon >0$ such that $u_{\tau-\epsilon} \geq u$ in $\R^N$. This will  contradict the definition of $\tau$ .
  
    By  the  estimate (\ref{eqholdloc})  in Theorem  \ref{IA},   there exists some constant $c>0$  such that for all $\epsilon >0$
              $$ |u_\tau-u_{\tau-\epsilon}|\leq  \epsilon   c. $$ 
Choosing $\epsilon$ small enough in order that $ \epsilon   c\leq {\eta\over 2}$ and $\epsilon < \tau$, one gets that  $u_{\tau-\epsilon}-u \geq 0$ on $\{x_1 = M_1\}$.           The same procedure as in  Claim 1 proves   that  the inequality holds in the whole space $\R^N$,  a contradiction with the definition of $\tau$. 
        \bigskip
        
        \noindent Hence
         $ \eta = 0$ and  there exists a sequence $(x_j)_j\in \left([-M_1, M_1] \times \R^{N-1}\right)^{\N}$ such that 
         $$(u-u_\tau) (x_j) \rightarrow 0.$$
Let $v_j (x) = u(x+ x_j)$ and $v_{j, \tau} (x) = u_\tau (x+ x_j)$; 
these are  sequences of bounded solutions, by uniform elliptic  estimates (consequence of Theorem \ref{IA}),   one can extract subsequences,  denoted in the same way,  such that  
$$v_j\rightarrow \bar v\quad\mbox{ and }\quad v_{j, \tau} \rightarrow \bar v_\tau$$ 
uniformly on every compact set of $\R^N$. 
Moreover,  $v_j$ and $v_{j,\tau}$ are solutions of the same equation and passing to the limit, $\bar v \geq \bar v_\tau$. 
Furthermore  $\bar v (0) = \lim_{j\rightarrow +\infty} u(x_j)= \lim_{j\rightarrow +\infty} u_\tau (x_j) = \bar v_\tau (0)$ and 
$$|\nabla  \bar v|(0) = \lim_{j\rightarrow +\infty}  |\nabla  u(x_j)| \geq m$$ 
by the assumption on $\nabla  u$.  
          
     Since $|\nabla  \bar v| >0$ everywhere,    by the  strong comparison principle  in Proposition \ref{strict},
  $\bar v_\tau = \bar v$ on any neighborhood of $0$ 
 . This would imply that $\bar v$ is $\tau$ periodic.   
 
 By our choice of $M_1$, 
 $\forall x\in \Sigma^+_{2M_1}$, 
 $v_j(x) = u(x+ x_j) \geq 1-\delta$ and 
 
   $\forall x\in \Sigma^-_{(-2M_1)}$, 
 $v_j(x) = u(x+ x_j) \leq -1+\delta$ , 
This contradicts the periodicity. Hence $\tau = 0$ and this ends the proof of Claim 2.

\bigskip
This implies that $\partial_{\vec \nu}  u(x) \geq 0$, for all $x \in \R^{N}$ since for all $t>0$, $u(x+t\vec \nu)\geq u(x)$ as long as $\nu_1>0$.  
 
 Take a sequence $\vec{\nu_n}=(\nu_{1,n}, \nu^\prime)$ such that $0<\nu_{1,n}$ and $\nu_{1,n}\rightarrow 0$.  Since $ u$ is ${\cal C}^1$,   by passing to the limit, 
 $$\partial_{\vec{\nu^\prime}} u(x) \geq 0.$$ 
 This is also true by changing $\vec \nu^\prime$ in $-\vec\nu^\prime$ , so finally $\partial_{\vec \nu^\prime} u(x) = 0$. 
 This ends the proof of Theorem \ref{th1}.

  \section{Existence's  results for the ODE.}
  We prove in this section that the one dimensional problem  (\ref{dim1wholespace}),  under additional assumptions on $f$, admits a solution and that,   when $\alpha \leq 0$, the solution is unique up to translation. 
 
 We consider the model Cauchy problem
    \begin{equation}\label{cauchydel}\left\{ \begin{array}{lc}
   -{\cal M}^+_{a,A} (u^{\prime\prime}) |u^\prime|^\alpha =  f(u), &  \mbox{in}\quad \R\\
                                u(0) = 0, u^\prime (0) = \delta  &
                                \end{array} \right.
  \end{equation}
where ${\cal M}^+_{a,A}$ is one of the Pucci operators.
  
With $f$ such that    
$f(-1) =f(0 ) = f(1)= 0$, $f$ is  positive in $]0 ,1[$, negative in $]-1, 0[$, $f$ is ${\cal C}^1([-1, 1])$.

\noindent We introduce the function $f_{a,A}(t)=\left\{\begin{array}{lc} \frac{f(t)}{a} & \mbox{if}\ f(t)>0\\
  \frac{f(t)}{A} & \mbox{if}\ f(t)<0
  \end{array}
  \right.
  ,$
so that equation (\ref{cauchydel}) can be written in the following way
   \begin{equation}\label{cauchydelta}\left\{ \begin{array}{lc}
   -u^{\prime\prime} |u^\prime|^\alpha =  f_{a,A}(u), &  \mbox{in}\quad\R\\
                                u(0) = 0, u^\prime (0) = \delta.  &
                                \end{array} \right.
  \end{equation}
  We  also assume on $f$:
\begin{enumerate}
\item $f^\prime(\pm1)<0$,    

\item $\displaystyle\int_{-1}^{1} f_{a,A}(s) ds =0$,                                               

\item for all $t\in (-1, 0]$,  $\int_t^{1 } f_{a,A}(s)ds >0$.
\end{enumerate}
$\delta_1$ will denote the positive real 
  \begin{equation}\label{delta1}
  \delta_1 = \left((2+\alpha ) \int_0 ^{1} {f(s)\over a} ds \right)^{1\over 2+\alpha}.
  \end{equation}
  Without loss of generality $f$ is extended outside of $[-1,1]$
 so that   $f\in {\cal C}^{0,1} (\R)$, $f\geq  0$ on $(-\infty, -1)$, $f\leq 0$ on $[1, +\infty)$. Then $ f$ satisfies also  for all $t\in \R\setminus\{\pm1\}$
$$\int_t^{1 } f_{a,A}(s) ds  >0.$$                                          
                                              
According to Cauchy-Lipschitz's theorem,  as soon as $u^\prime (0) \neq 0$ there exists a local unique solution.  Moreover the Cauchy Peano's Theorem establishes some   global existence's theorem.

We establish existence and uniqueness  (in the case $\alpha \leq 0$) of  weak solutions and their equivalence with  viscosity solutions.  
\begin{defi}
  A weak solution for (\ref{cauchydelta}) is a ${\cal C}^1$ function which satisfies  in the distribution sense 
               \begin{equation}\label{weak}\left\{ \begin{array}{lc}
               -{d\over dx} (|u^\prime|^\alpha u^\prime) = (1+\alpha) f_{a,A}(u)&\mbox{ in }  \ \R\\
               u(\theta) =0, \ u^\prime(\theta) = \delta.& \ 
               \end{array}\right.
                \end{equation}
Without loss of generality we can suppose that $\theta=0$. \end{defi}
Remark that we are interested in solutions that are in $[-1,1]$ so we shall suppose that $u_o\in (-1,1)$.

\begin{rema} 
 Let us note that the condition 2 on $f$  is necessary for the existence of  weak solutions  which satisfy $\lim_{x\rightarrow +\infty} u(x) = 1$, $\lim_{x\rightarrow -\infty} u(x) = -1$. Indeed    by continuity $u$ has a zero and without loss of generality we can suppose that it is in 0.  Since the solution $u$ is ${\cal C}^1$, and  bounded, the limit of $u^\prime$ at infinity is $0$. 
 In particular, multiplying the equation (\ref{weak}) by $u^\prime$ and integrating in $[0,+\infty)$  
  $$|u^\prime (0 )|^{2+\alpha}=  -(2+\alpha) \int_0^{1} {f(s)\over a} ds$$
  and  in $]-\infty,0]$, 
  similarly  
  $$ |u^\prime (0 )|^{2+\alpha} =(2+\alpha) \int_{0}^{-1} {-f(s)\over A} ds= (2+\alpha) \int_{-1}^0 {f(s)\over A} ds .$$
This implies 2.
\end{rema}

\begin{prop} For $\alpha>-1$ there exists a solution of (\ref{weak}), and for  $\alpha\leq 0$ this solution is unique.\end{prop}     
Proof. 

To prove  existence  and uniqueness observe that   both the equations (\ref{cauchydelta}) and (\ref{weak}) can be written, with $u=X$ and $Y= |u^\prime|^\alpha u^\prime$, under the following  form 
                \begin{equation}\label{eqcauhlip} \left(\begin{array}{c} X^\prime\\
                Y^\prime\end{array} \right) =    \left( \begin{array}{c}
                                                   |Y|^{\frac{1}{\alpha+1}-1}Y\\
                                    -(1+\alpha)f_{a,A}(X)                                     \end{array}\right)
                                    \end{equation}
with the initial conditions $X(0) = 0$, $Y(0) =| \delta|^\alpha \delta $ and the map $(X, Y) \mapsto  \left( \begin{array}{c}
                                                      |Y|^{\frac{1}{\alpha+1}-1}Y\\
                                    -(1+\alpha)f_{a,A}(X) 
                                    \end{array}\right)$
is continuous. When $\alpha\leq 0$  it is Lipschitz continuous; and when $\alpha>0$ it is  Lipschitz continuous for $Y(0)\neq 0$.
Now the result is just an application of the classical Cauchy Peano's  Theorem,
and the Cauchy Lipschitz theorem.
It is immediate to see that weak solutions and the  solutions  of (\ref{eqcauhlip}) are the same.
This ends the proof.

\bigskip
      Observe that weak solutions are viscosity solutions. Indeed,  it is clear that $|u^\prime |^\alpha u^\prime$ is ${\cal C}^1$, hence if $u^\prime \neq 0$,  $u^\prime$ is ${\cal C}^1$. Finally $u$ is ${\cal C}^2$ on each point where the derivative is different from  zero  and on such a point the equation is  $-|u^\prime |^\alpha u^{\prime\prime} = f(u(x))$ so $u$ is a viscosity solution. 
      
      We now consider the case where $u$ is locally constant  on $]x_1-\delta_1, x_1+\delta_1[$ for some $\delta_1  >0$ the "weak  equation" gives $f(u(x_1))= 0$, then $u(x_1) = 0, 1$ or $-1$, and $u$ is a viscosity solution.

     We  now assume that $\alpha \leq 0$ and  recall that according to the regularity results in \cite{BD10} applied in the one dimensional case, the solutions are ${\cal C}^2$. We now prove that the viscosity solutions are weak solutions.   
 
When $u^\prime(x)\neq 0$ or when $u$ is locally constant,  it is immediate that $u$ is a weak solution in a neighborhood of that point. 

So, without loss of generality, we suppose that,  
$u^\prime (x_1) = 0$,  $1>u(x_1) >0$ and hence $u$ is not locally constant.
Then, by continuity of $u$ and the equation, there exists $r>0$ such that
$$u^{\prime\prime}\leq 0 \quad\mbox{in}\quad (x_1-r,x_1+r).$$
Furthermore there exists $(x_n)_n$, such that $x_n\in (x_1-r,x_1)$, $x_n\rightarrow x_1$ and 
$u^\prime(x_n)\neq 0$;
by the equation we obtain that 
$$u^{\prime\prime}(x_n)<0.$$
Finally, $u^\prime(x)=\int_{x_1}^x u^{\prime\prime}(t)dt>0$ for $x\in (x_1-r,x_1)$.
Similarly $u^\prime(x)<0$ for $x\in (x_1,x_1+r)$.

By uniqueness of the weak solutions , $u$ satisfies in a neighborhood of $x_1$:
  $$-{d\over dx} (|u^\prime|^\alpha u^\prime) = {(1+\alpha)f(u(x))\over a}.$$    
This proves that $u$ is a weak solution. 
 
 \begin{prop}\label{propcauchy}  
Suppose that $\alpha \leq 0$. Let $u_\delta$ be the unique solution  of  (\ref{cauchydel}). Then for  $\delta_1$ defined in (\ref{delta1}),      
      
\noindent 1) If $\delta > \delta_1$,     $|u_\delta (x)| \geq C|x|$ for $C=\delta^{2+\alpha}-\delta_1^{2+\alpha}$.
In particular $\displaystyle\lim_{x\rightarrow  \pm \infty}u_\delta(x)  = \pm \infty$   and $u^\prime_\delta >0$.                               

\noindent 2) If $\delta = \delta_1$,  $u^\prime _\delta >0$ in $\R$ and $\displaystyle \lim_{x\rightarrow +\infty} u_\delta (x)  = 1 $, $\displaystyle\lim_{x\rightarrow  -\infty}u_\delta (x)  = -1$.                                  

\noindent3) If $-\delta_1\leq \delta < \delta_1$ then $ |u_\delta(x)|_\infty < 1$ for any  $x\in\R$.  The solution can oscillate. 
                                 
\noindent 4) If $\delta < -\delta_1$, $u_\delta $ is decreasing on $\R$, hence $u_\delta<0$ on $\R^+$, $u_\delta >0$ on  $\R^-$.

\end{prop}

\begin{rema}
The  case 2) in  Proposition  \ref{propcauchy} is clearly false in the case $\alpha >0$.  As one can see with the example : 
$\alpha = 2$, $f(u) = u-u^3$, $u(x) = \sin x$, $u
$ satisfies $u^\prime (0) = \delta_1 = 4 \int_0^1 f(s) ds$, $u(1)={\pi\over 2}$ and  it oscillates.  

\noindent However the   conclusion in the other cases    holds   for any $\alpha$.
\end{rema}
Proof of Proposition \ref{propcauchy}.

1 \& 4)  To fix the ideas we suppose that $\delta> \delta_1$, the proof is identical in the case $\delta<-\delta_1$. For   $x>0$, since $u_\delta>0$   one has 
 \begin{eqnarray*} |u^\prime_\delta| ^{2+\alpha} (x) &=& \delta^{2+\alpha} - (2+\alpha) \int_0^{u_\delta(x)} {f(s)\over a} ds \\
 &=& \delta^{2+\alpha} -\delta_1^{2+\alpha}    + (2+\alpha)\int_{u_\delta (x)}^{1} {f(s)\over a} ds \\
&\geq & \delta^{2+\alpha} -\delta_1^{2+\alpha}:=C.
\end{eqnarray*}
This proves, in particular, that $u_\delta ^\prime (x) \neq 0$ for all $x$ and the Cauchy Lipschitz theorem ensures the local existence and uniqueness on every point, hence also the global existence . From  this, we also derive that  $u^\prime _\delta >0$ and  for $x>0$,  $u_\delta (x ) \geq Cx$, and symmetric estimates for $x<0$ give $u_\delta (x) \leq  C x$.                                   
                                   
 2)   If $\delta = \delta_1$ then   $|u^\prime_\delta|^{2+\alpha}(x)  = (2+\alpha )\int_{u_\delta (x)} ^{1} {f(s )\over a} ds > 0$. Suppose that there exists some point $\bar x$ such that $u_\delta(\bar x) = 1$ then $u_\delta^\prime(\bar x ) = 0$.
By  the uniqueness of the solution  $u_\delta(x)\equiv 1$ which contradicts the fact that $u_\delta^\prime (0) = \delta_1\neq  0$.  
  
   We have obtained that $u_\delta (x) < 1$ everywhere. Moreover $u_\delta $ is increasing and bounded  then $\lim_{x\rightarrow+  \infty}  u_\delta^\prime = 0$. By hypothesis 3.  on  $f$, this  implies that $\lim_{x\rightarrow  + \infty}  u_\delta (x) =   1$.

3)   Suppose that $0 < \delta < \delta_1$,    and let $\theta^+$  be  such that $(2+\alpha) \int_0^{\theta^+}  \frac{f(x)}{a}dx = \delta^{2+\alpha}$, which exists by  the mean value theorem.  Either $u_\delta < \theta^+ $ for all $x$, or there exists $x_1$ such that $u_\delta (x_1) = \theta^+$, and then $u^\prime_\delta (x_1) = 0$.  Let us note that                                      $ u = \theta^+ $ on a neighborhood of $x_1$ is not a solution since  $f(\theta^+) \neq 0$. 
So $u_\delta$ is not locally constant and in particular, in a right neighborhood of $x_1$:

$$ \exists \varepsilon_o, \ u^{\prime\prime}_\delta (x) \leq 0, \ u^{\prime\prime}_\delta\not\equiv 0$$
for all $x\in (x_1, x_1+\varepsilon_o)$, hence 
$u^\prime_\delta (x) <0$ in $(x_1, x_1+\varepsilon_o)$. 

So $u$ is decreasing until it reaches a point where $u^\prime_\delta (x_2) = 0$. Observe that by the equation 
$$ 0 =  |u^\prime_\delta| ^{2+\alpha} (x_2) = - (2+\alpha) \int_{\theta^+}^{u_\delta(x_2)} f_{a,A}(s) ds .$$
Hence $u(x_2) = \theta^-\in (-1,0)$. 

We can reason as above  and obtain that $u$ oscillates between $\theta^-$ and $\theta^+$.

\end{document}